\documentclass{article}
\usepackage[utf8]{inputenc}

\title{Non cyclic division algebras of prime degree}
\author{Shmuel Rosset\\
Tel Aviv University }
\date{September 9, 2020}

\usepackage{graphicx}

\begin{document}

\maketitle

\section{Introduction}

If $D$ is a division algebra of degree $n$, i.e.\ dimension $n^2$ over its center $k$, then it 
contains maximal commutative subfields which are separable extensions of $k$ of dimension $n$. 
If there is such a maximal subfield, $K$, which is a Galois extension of $k$ we say that $D$ is a 
{\em crossed product}. The Skolem-Noether
theorem says that every element of the Galois group $G=G(K/k)$ can be extended to an inner 
automorphism of $D$. The Galois group becomes a kind of "Weyl group" in the sense that it is 
$N_{D^*}(K^*)/K^*$, where $N_{B}(A)$ denotes the normalizer of $A$ in $B$. 
The group extension $$1\to K^*\to N_{D^*}(K^*)\to G(K/k)\to 1$$
determines a class in $H^2(G,K^*)$ and it also determines $D$ and its Brauer class $[D]$ in Br$(K/k)$,
the subgroup of Br$(k)$ of the elements split by $K$.

A crossed product in which the Galois group is cyclic is called a cyclic algebra. Following Hamilton's
quaternions the first division algebras were cyclic algebras. Remarkably it turned out, proved by 
Merkurjev and Suslin \cite{MS}, that in the presence of roots of unity cyclic algebras generate the 
Brauer group of a field.
Still the question was asked: is every 
division algebra a crossed product? In other words, does every 
division algebra contain a maximal subfield which is a Galois extension of the center? The first to 
construct division algebras that are not crossed products was Amitsur \cite{Am} who showed that his generic matrix
algebras of degree $n$ are not crossed products when $n$ is divisible by the square of an odd prime or by 8.
There have been other constructions since but none of non-cyclic algebras of prime degree. Here this is done, by a completely different method. The existing examples of non-crossed products are all of explicitly constructed
algebras, and the proof that they are not crossed product requires work. In this paper the algebras that are 
not cyclic are subalgebras, over the same center, of certain crossed products. 
These algebras are known
to exist but no explicit description for them seems to be known. 
Yet to prove that they are not cyclic, or crossed products, is elementary. 

The subalgebras we are talking about are
the primary components of division algebras that are crossed products of groups that have 
few, or hardly any, homomorphic 
images such as simple non-cyclic groups or full symmetric groups. 

Given a finite group $G$ of order $n$,
the existence of division algebra crossed products of degree $n$ with Galois group $G$, in any 
characteristic and free of any assumption on roots of unity, appears as 
the "generic" crossed products that were constructed
many years ago in \cite{Ro}. In fact, if $m$ divides $n$ and they have the same radical,
i.e.\ the same set of prime divisors,
then a division algebra crossed product with group $G$, of degree $n$ and order $m$ in the Brauer group 
of its center, is constructed there. In this paper the order of the division algebras will not play a part,
but it is perhaps noteworthy that for smaller $m$ centers of larger dimension are required.

In \cite{Ro} the building materials for constructing a generic
$G$ crossed product are taken from a free
presentation of $G$:  $$1\to R\to F\to G\to 1,$$ 
the main building block being the relation module $R_{ab}:=R/[R,R].$
It seems that relation modules of finite groups were of interest from the early days of the cohomology of groups.
See the original paper
of Eilenberg and Maclane \cite{EM} where relation modules occupy center stage.

The presentation gives rise to an
extension $ F/[R,R]\to G\to 1$ whose kernel, the relation module
$R_{ab}$, is a free Abelian group and a $G$ lattice. 
The group $F/[R,R]$ is always torsion free and its group ring, over a base field $\ell$, 
is a Noetherian {\em domain}  whose 
classical ring of fractions is a finite dimensional division ring. This division ring 
is what we call a generic crossed product
with group $G$. The action of $G$ on $A=R_{ab}$ is faithful (assuming $R$ is not cyclic) and 
induces an action on the field $\ell(A)$, the field of 
fractions of the group ring $\ell A$. Thus the center of this division algebra is the fixed field $\ell(A)^G.$

In this paper we show that non-crossed products, and even non-cyclic algebras of prime degree,
can hide in plain sight as primary components of such generic crossed products. 
Primary components of crossed products, in particular crossed products 
that arise from localising "prime" group rings
(defined in \S3) of virtually free abelian groups, are 
discussed in \S3. Though these primary components are perhaps mysterious they do have one very 
important property. Suppose your division algebra, $D$, is a
crossed product of the Galois extension $K/k$ whose Galois 
group is $G$, $P$ is a $p$-Sylow subgroup of $G$ and $D(p)$ is the $p$ primary component of $D$. 
This primary component is a division algebra with center $k$ of degree $|P|$. 
It exists in some, possibly high, power of $D$, 
but we don't see it. 
However its "restriction" to $P$, one manifestation of which takes the class of $D(p)$, in Br$(k)$,  
to the class of $K^P\otimes_k D(p)$ in Br($K^P$), 
is the class of the division algebra crossed product
of $K/K^P$ and $P$, with factor set the restriction to $P$ of the factor set utilized to obtain $D$. 

It is this property that enables us to show that in a
division algebra which is a crossed product with a "difficult" (i.e.\ with few non-trivial homomorphic images) 
Galois group, the 
primary components $D(p)$ cannot
be crossed products, even when the Sylow group is cyclic of order $p$.
In fact when the Sylow subgroup is cyclic the assumptions on the Galois group are weaker than the assumptions 
needed in the general case. So we separate the two cases.\\
{\bf Theorem 1.} {\em Let $p$ be an odd prime, and $G$ a finite group whose $p$ Sylow subgroup is of order 
$p$ but which does not have a normal subgroup of index $p$. If $D$ is a division algebra crossed product of a 
Galois extension $K/k$ with Galois group $G$ and an appropriate factor set
then the $p$ primary component of $D$, which is a division algebra of degree $p$ central over $k$, is not cyclic.}

A more general statement, applying to all odd Sylow subgroups, is true but for a smaller category of groups.\\ 
{\bf Theorem 2} {\em  Suppose $G$ is a finite group that has a non-commutative simple subgroup of index $\leq 2$
and $D$ is a division algebra crossed product with Galois extension $K/k$, Galois group $G$ and suitable
factor set. If $p$ is an odd prime such that $v_p(|G|)=a$ for some $a\geq 1$, then
the $p$ primary component of $D$, 
which is a division algebra of degree $p^a$
over $k$, is not a crossed product.} 

The case $p=2$ in theorem 2 is left open for now.

Throughout this paper the expressions "cocycle" and "factor set" will both be used for the same thing. Cocyles
are for group extensions what factor sets are for central simple algebra crossed products.

\section{Relation modules}

As noted in the introduction
if $F$ is a free group and $$1\to R\to F\to G\to 1$$ is a free presentation of the group $G$ then 
$R_{ab}=R/[R,R]$ is a $G$ module, which is called the relation module.  
We are interested in the group $F/R'$ when $G$ is finite. Presentations are far from unique, of course,
so every presentation carries its own relation module, but they all have the same
cohomology. In fact applying Tietze transformations shows that all relation modules of a finite group
are stably isomorphic:
if $M,N$ are two relation modules for $G$ there are
finitely generated free ${\bf Z}G$ modules $E,F$ such that $M\oplus E$ and $N\oplus F$ are isomorphic.
Since free modules are cohomologically trivial we see that the Tate cohomology of the relation module is
uniquely determined. It is also easy to see now that, unless $F$ is cyclic, the action of $G$ on the relation 
module is faithful. To show that
suppose $1\neq x\in G$ and let $C_x$ be the cyclic subgroup 
generated by $x$. The inverse image of $C_x$ in $F$ gives a presentation with kernel $R$ and thus 
the relation module is $R_{ab}$. 
But as a relation module for $C_x$, i.e.\ as a ${\bf Z}C_x$ module, it is a direct
sum of a trivial module, coming from the presentation ${\bf Z}\to C_x\to 0$, and a positive number of
free modules. And clearly $x$ acts non-trivially on a free $C_x$ module. The fact that $F/R'$ is torsion 
free can also be proved at this point but is given a somewhat different proof below. 

A similar proof shows that if $G$ is not cyclic then the center of $F/R'$ is trivial. Indeed, let $x,y$
be two elements in $G$ such that the group they generate is not cyclic. The invariants of $R_{ab}$
under the action of $C_x$ are elements of the infinite cyclic group generated by $x$ in $F/R'$. These are not invariant under $y$, which proves the claim. 

Let $\Delta(G)$ denote the augmentation ideal so that 
$$0\to \Delta(G)\to {\bf Z}G \to{\bf Z}\to 0$$ is an exact sequence of $G$ lattices. If $G$ can be generated 
by $d>1$ elements we can take $F$ to be free of rank $d$. Then there is, less obviously, 
an exact sequence of $G$ modules $$0\to R_{ab}\to {{\bf Z}G}^{\oplus d} \to\Delta(G)\to 0$$
where ${{\bf Z}G}^{{\oplus}d}$ denotes a free ${\bf Z}G$ module of rank $d$. See \cite{DLJ} Ch. 11.

With ${\hat H}(~,~)$ denoting Tate cohomology,
the first exact sequence implies a natural isomorphism 
${\hat H}^n(G,{\bf Z}) \cong {\hat H}^{n+1}(G,\Delta(G)),$ while the second implies an isomorphism
${\hat H}^n(G,\Delta(G)) \cong {\hat H}^{n+1}(G, R_{ab})$ which is also,
but again less obviously, natural. We only need
the isomorphism $${\hat H}^2(G,R_{ab}) \cong {\hat H}^0(G,{\bf Z})\cong {\bf Z}/|G|{\bf Z}.$$
It is shown in \cite{Ro} that the extension  
$$ \alpha: 1\to R_{ab}\to F/R'\to G\to 1,$$ derived from the given presentation, generates the cyclic group 
${\hat H}^2(G,R_{ab}).$ Indeed let $\beta: 1\to R_{ab}\to E\to G\to 1$ be a generator of
${\hat H}^2(G,R_{ab}).$ Using the freeness of $F$ one shows that
there is a map $f: F/R' \to E$ such that $f_{*}(\alpha)=\beta.$ Thus the order of $\alpha$ is $|G|$ and it is 
also a generator.

If $H$ is a subgroup
of $G$ its inverse image, $F_H$, in $F$ is a presentation of it with kernel $R$ and, by the same token,
${\hat H}^2(H,R_{ab})$ is cyclic of order $|H|$ and the extension $1\to R_{ab}\to F_H/R'\to H\to 1$
is a generator. This proves that $F/R'$ is torsion free because if it had torsion then for some cyclic
subgroup $H$ the extension would be split, which we know is not the case.

The group ring of $F/R'$, over a base field $\ell$, contains the group ring of the free Abelian group $R_{ab}$.
The action of $G$ on $R_{ab}$ extends to an action on $\ell R_{ab}$ and on its field of fractions
$\ell(R_{ab})$. Denote $\ell(R_{ab})$ by $K$ and its fixed subfield under the action of $G$, $K^G$, by $k$.

It is easy to see, and proved in \cite{Ro}, that $R_{ab}$ is a direct summand (as ${\bf Z}G$ modules)
of $K^*$ and hence the inclusion $\iota: R_{ab}\hookrightarrow K^*$ induces an injection on the cohomology. It 
follows that the cohomology class $\iota_{\ast}(\alpha)$ is of order $n=|G|$.
The crossed product of $K/k$ and $G$ with the cocycle defining $\alpha$ is therefore a central simple 
algebra of degree $n$ (i.e.\ dimension $n^2$ over $k$) whose Brauer class is of order $n$. Hence it is a 
division algebra. And it is obviously also the "classical" division ring of fractions of the group ring
$\ell[F/R']$. Note that 
our knowledge that $\alpha$ is of order $n$ implied that the group ring is a domain.

In fact if $\Gamma$ is a virtually abelian torsion free group then $\ell\Gamma$ is a domain for every 
field $\ell$. This is a non-trivial statement first proved in \cite{FS}, for $\ell$ of characteristic 0, 
and in general in \cite{Li}.
The total classical ring of fractions exists and is a division ring. If, moreover, $\Gamma$ has a normal 
commutative subgroup of finite index $C$ such the action of $\Gamma/C$ on $C$ is faithful then
the division ring of fractions is the crossed product of the field $\ell(C)$ and $\Gamma/C$ with the cocycle
provided by the extension $1\to C\to\Gamma\to \Gamma/C \to 1$. As we knew that $F/R'$ is torsion free
we see that the information that the order of the extension $\alpha$ is $n=|G|$ is actually redundant.

\section{primary components of division algebras}

Let $D$ be a division algebra over the field $k$ of degree $n$ and suppose $n=rs$ where $r,s$ are relatively
prime and both greater than 1. 
If $a,b$ are integers such that $ar+bs=1$ modulo $n$, 
what can be said of $D^{ar}$ i.e.\
$D\otimes_k\cdots \otimes_{k} D$ $ar$ times? 
Wedderburn's theorem tells us that it is isomorphic to some $M_{\nu}(D')$
with $D'$ a division algebra over $k$. While $D'$ is uniquely determined there is little we can say on $\nu$ . 
Similarly $D^{bs}\approx M_{\mu}(D'').$
What are the degrees of $D', D''$? Since $ar=1$ modulo $s$ and the Brauer
class $[D]=[D^{ar}] [D^{bs}]$ 
the best we can hope for is $s$ for $D^{ar}$ and $r$ for $D^{bs}$. And 
indeed this is precisely the case; see, for example, \cite{Sa} ch. 5 where an elaborate proof using symmetrizers 
is worked out.

It follows from these considerations that if the division algebra
$D$ is a crossed product of $K/k$ with group $G$ of order $n=p_1^{e_1}\cdots p_v^{e_{v}}$, where $p_i$ are 
primes, then there are division algebras over $k$, $D(p_i)$ of degree 
$p_i ^{e_i},~i=1,...,v$ such that $D\approx \otimes_{i=1}^v D(p_i)$. Explicitly, let $n_i= n/p_i ^{e_i}$
and $a_i$ integers such that $\sum_{i=1}^v a_i n_i=1$. Then $D(p_i)$ is the division algebra component of
$D^{a_i n_i}$. These are the primary components of the 
division algebra.

Let $P_i$ denote a Sylow $p_i$ subgroup of $G$. 
The cohomology class  corresponding to
$D(p_i)$ being a power of the cohomology class corresponding to $D$, its restriction to $P_j$ is a power
of that of $D$. If $j\neq i$ then this restriction is zero as the exponent $a_i n_i$ is divisible by $p_j ^{e_j}$.
On the other hand, as $a_i n_i=1$ modulo $p_i ^{e_i}$ the restriction to $P_i$ is the "identity", i.e.\ 
it is represented by the division algebra crossed product of $K/K^{P_i}$ with Galois group $P_i$ and factor set
which is the restriction from $G$. In other words, the restriction 
of the $p$ primary component to the $p$ Sylow subgroup is the same as the restriction from $G$ to the Sylow subgroup.

Much more can be proved
when dealing with division algebras, and even central simple algebras, that are total classical rings of 
fractions of "prime" group rings of virtually abelian groups. For our purpose the relevant groups are
groups that are sometimes called "crystallographic":
extensions of a finite group $G$ by a torsion free ${\bf Z}G$ module $M$, fitting into an exact sequence 
$$1\to M\to E\to G\to 1$$  
such that the action of $G$ on $M$, 
by conjugation in $E$, is a faithful representation of $G$. 
This is the same as saying that $M$ is a maximal abelian
subgroup of $E$ and that $E$ is {\em prime} in the sense that it has no non-trivial finite normal subgroup.
In this case the group ring, over an integral domain, is a prime ring and
the total classical ring of fractions of the group ring $R:=\ell E$
is a simple artinian ring which is finite
dimensional over its center, i.e.\ a central simple algebra. Here $\ell$ is any field.
The center, as before, is the fixed field under 
the action of $G$ on the field of fractions $\ell(M)$. Wedderburn's theorem tells us that
the total ring of fractions, denoted $\ell(E)$, is a matrix algebra $M_{\nu}(D)$ where $D$ is a division algebra
with center $k=\ell(M)^G$. The number $\nu$ is the "Goldie rank" of the group ring $R$. 
For example, if the extension splits, i.e.\ the cohomology class associated to it is zero, then $D=k$ and 
$\ell(E) \approx M_n(k)$ where $n=|G|$. 
On the other extreme if $E$ is torsion free then, by the theorem of Farkas-Snider-Linnell quoted above,
$\ell(E)$ is a division ring and $\nu=1$. 

In \cite{R1} it was shown that the Goldie rank is equal to another number associated with the group ring $R$.
If $T$ is a finitely generated
$R$ module it is also finitely generated over the subring $\ell M$, which is a commutative Laurent 
polynomial ring. As polynomial rings are smooth $T$ has a finite projective resolution, i.e.\ an exact sequence
of $\ell M$ modules 
$$ 0\to Q_r \to \cdots \to Q_1\to Q_0 \to T\to 0 $$
in which the $Q_i$ are finitely generated projective $\ell M$ modules. 
The rank of an $\ell M$ module $Q$ is defined
to be the dimension over the field of fractions $\ell(M)$ of $\ell(M)\otimes_{\ell M} Q$. The Euler
characteristic of $T$, denoted $\chi_R(T)$,
is, by definition, $$\frac{1}{|G|}\sum_{i=0}^r (-1)^i {\rm rank}_{\ell M}(Q_i).$$
It is independent of the resolution. And it turns out that it is independent of 
the field $\ell$. In fact it depends only on the group $E$ and not on the extension, in the sense that any
subgroup $M'$ of finite index of $E$ such that $\ell M'$ is equally "smooth" would give the same result.
It is shown in \cite{R1} that the Goldie rank is equal to
the smallest positive integer $\rho$ such that $\rho\cdot \chi_R(T)$ is an integer for all 
finitely generated $R$ modules $T$. 

It is easy to see that if $H$ is a finite subgroup of $E$ then $$\chi_R(\ell[E/H])=\frac{1}{|H|},$$ where
$\ell[E/H]$ denotes the permutation module on the cosets of $H$. (If $H$ is an infinite subgroup 
$\chi_R(\ell[E/H])=0$.) Since the finite subgroups of
$E$ are just the subgroups of $G$ over which the extension splits the least common multiple
of orders of the finite subgroups divides $\rho$ and in \cite{R1} it was conjectured that this is an equality,
i.e.\ that $\rho$ is the least common multiple of the orders of finite subgroups. This was proved by Moody
\cite{Mo} who proved the stronger result that the permutation modules generate $G_0(R)$.

Moody's result makes possible the evaluation of the degree of the division algebra component of 
central simple algebras that are  
classical rings of fractions of crystallographic group rings.

\section{Proofs} 
This paper grew from the question 
"If $p$ is a prime, is the $p$ primary component of the generic crossed product division 
algebra, with Galois group the full symmetric group ${\mathcal S}(p)$, cyclic?" 
Theorem 1 is the negative answer.\\
{\bf Proof}:
Let $p$ be a prime $\geq 3$ and $G$ a group whose $p$ Sylow subgroup is cyclic of order $p$ but it does not
have a normal subgroup of index $p$. The symmetric groups ${\mathcal S}(p)$ satisfy these assumptions, even
${\mathcal S}(3)$. Let $D$ be a division algebra which is a crossed product of a Galois extension $K/k$
whose Galois group is $G$, with an appropriate factor set. 

The examples, in \S2, for such crossed products are the rings of fractions of group rings of 
torsion free groups $E$ that are extensions 
$$ 1\to A\to E\to G\to 1$$ in which $A$ is abelian and is faithful as a $G$ module. 
The generic such extensions are
those that come from free presentations of $G$, the module $A$ being the relation module associated with the
presentation.

Let $H$ be a
$p$ Sylow subgroup of $G$. By assumption it is cyclic of order $p$. 
Its fixed field is $K^H$ and the crossed product of $K/K^H$ and $H$, the factor set being the restriction to $H$
of the factor set defining $D$,
is a cyclic algebra of 
degree $p$ (over its center $K^H$) within $D$. We denote it $D_H$. 

Since primary components are uniquely determined up to isomorphism we will refer to them as known. So
let $D(p)$ be the $p$ primary component of $D$ over $k$. It is a division algebra
of degree $p$ with center $k$ such that the Brauer class of $D$ is the product of the Brauer 
class of $D(p)$ and another class of order prime to $p$. Being a crossed product
the cohomology class representing $[D]$ is in $H^2(G,K^*)$. The Brauer class of $D(p)$ is a 
power of that of $D$ and as such it is also in $H^2(G,K^*)$. Thus we can take its restriction to $H^2(H,K^*)$.
What we know is that this restriction 
to $H$, or equivalently $[K^H\otimes_k D(p)]$, is equal to $[D_H]$. 

We will show that this cannot hold if $D(p)$ is a crossed product, i.e.\ cyclic in this case.
Suppose, by contradiction, that $D(p)$ 
is a cyclic algebra. This means that there is a cyclic extension of degree $p$, $L/k$,
in $D(p)$ which makes it into a cyclic division algebra. Now $K$ and $L$ are two Galois extensions of $k$
both subfields of a given separable closure, $k_s$, of $k$. If $Y$ is a finite extension of $k$ 
contained in $k_s$ we denote the Galois group of $k_s/Y$ by ${\mathcal G}_Y$.
In this notation the absolute Galois group of $k$ is ${\mathcal G}_k$, and it
has two normal subgroups of finite index: ${\mathcal G}_K$ and  ${\mathcal G}_L$ with quotients 
identifiable, via the restriction of Galois action map, with the Galois groups
$G(K/k),~G(L/k)$ respectively.
Now ${\mathcal G}_L$ cannot contain ${\mathcal G}_K$ because if it did then $G(K/k)$ would have a 
normal subgroup of index $p$, which, by assumption, is not the case.

It follows that ${\mathcal G}_{KL}$, which is equal to 
${\mathcal G}_K \cap{\mathcal G}_L$, is a proper subgroup of ${\mathcal G}_K$.  
In fact 
$$({\mathcal G}_k :{\mathcal G}_K \cap{\mathcal G}_L) = ({\mathcal G}_k :{\mathcal G}_K)\cdot
({\mathcal G}_k :{\mathcal G}_L).$$
To prove that divide by ${\mathcal G}_{KL}$. The group ${\mathcal G}_k/{\mathcal G}_{KL}$ is the Galois group
$G(KL/k)$ and it has 2 normal subgroups ${\mathcal G}_K /{\mathcal G}_{KL}$, identifiable as $G(KL/K)$, and 
${\mathcal G}_L /{\mathcal G}_{KL}$ identifiable as $G(KL/L)$. The intersection, $G(KL/K)\cap G(KL/L)$ in 
$G(KL/k)$ is trivial and it follows that these subgroups commute elementwise, i.e.\ every element in one commutes 
with every element in the other. Let $\pi_1: G(KL/k)\to G(L/k)$ be the restriction map. It is surjective and 
its kernel is $G(KL/L)$. The restriction of $\pi_1$ to $G(KL/K)$ is an injection to $G(L/k)$ with
a non-trivial image. Since $G(L/k)$ is of order $p$ it is surjective and the restriction of $\pi_1$ to 
$G(KL/K)$ is an isomorphism $G(KL/K)\approx G(L/k)$.
Similarly, if $\pi_2: G(KL/k)\to G(K/k)$ is the restriction map it induces an isomorphism 
$G(KL/L) \approx G(K/k)$.

The map $\pi: G(LK/k)\to G(L/k)\times G(K/k)$ defined by
$\pi(x)=(\pi_1(x),\pi_2(x))$ is
obviously injective and as $|G(KL/k)|=|G(K/k)|\cdot|G(L/k)|$ it is an isomorphism. It identifies $G(K/k)\times 1$
as the kernel of $\pi_1$. 
The map induced by $\pi_1$ in cohomology is the {\em inflation} map
$$ {\pi_1}^* : H^2(G(L/k),L^*)\to H^2(G(KL/k),(KL)^*).$$ It is injective and its image is equal to the kernel
of the restriction map $$H^2(G(KL/k),(KL)^*\to H^2(G(K/k),K^*).$$ This is the well known inflation-restriction 
exact sequence, see \cite{GS} p.88. 

Thus the restriction to $G(K/k)$ of every element which is an inflation from $H^2(G(L/k), L^*)$ is zero. In particular the same is true for restriction to $H$ which is a subgroup of $G(K/k)$. It follows that if
$D(p)$ is cyclic it cannot restrict to a non-trivial element in $H^2(H,K^*)$, as it must. This shows that $D(p)$
is not a cyclic algebra, proving theorem 1.

The proof of theorem 2 is similar. We use the notation of the theorem. Let $H$ be a $p$ Sylow subgroup,
$D_H$ the division subalgebra of $D$ with center $K^H$ which is the crossed product of $K/K^H$ and $H$, with
factor set the restriction from $G(K/k)$. As before $D(p)$ will denote the $p$ primary component of $D$. 

If $D(p)$ is a crossed product it has a maximal commutative subfield $L$, of dimension $p^a=|H|$ over $k$, 
which is a Galois extension of $k$. As before if $Y\subset k_s$ is a finite extension of $k$ we denote the
Galois group of $k_s/Y$ by ${\mathcal G}_Y$. Then $${\mathcal G}_K\cap {\mathcal G}_L={\mathcal G}_{KL}.$$
The compositum $KL$ is not equal to $K$ because if it were then ${\mathcal G}_K \subset {\mathcal G}_L$ which implies that $G=G(K/k)$, which we identify with ${\mathcal G}_k/{\mathcal G}_K$, has a normal subgroup 
${\mathcal G}_L/{\mathcal G}_K$. Our assumption was that the only non-trivial normal subgroup of $G$ is of
index $\leq 2$. The index of ${\mathcal G}_L$ in ${\mathcal G}_k$ is the order of a Sylow $p$ subgroup
of $G$ which is certainly not 2. It follows that $KL\neq K$, as claimed.

Thus $G(KL/k)$ contains two non-trivial normal subgroups, $G(KL/K)$ and $G(KL/L)$ whose intersection is trivial,
which implies that they commute elementwise. If $\pi_1: G(KL/k)\to G(K/k)$ is the restriction map, with kernel
$G(KL/K)$, and $\pi_2: G(KL/k)\to G(L/k)$ is the restriction map with kernel $G(KL/L)$, then the map
$$\pi: G(KL/k)\to G(K/k)\times G(L/k),~~\pi(x)=(\pi_1(x), \pi_2(x)),$$ is injective. 

We will prove that when $p$ is odd $\pi$ is an isomorphism. 
The restriction of $\pi_1$ to $G(KL/L)$ is injective since its
intersection with the kernel of $\pi_1$ is trivial. The image of $G(KL/L)$ in $G(K/k)$ is a non-trivial
normal subgroup and hence either the whole of $G(K/k)$ or a subgroup of index 2. Suppose it is of index 2. 
The equality $|G(KL/K)|\cdot|G(K/k)|=|G(KL/L)|\cdot|G(L/k)|$ leads to the conclusion that $2|G(KL/K)|=|G(L/k)|$.
But $G(L/k)$ is a $p$ group and $p$ is odd, which is impossible. Thus $\pi_1$ induces an isomorphism of 
$G(KL/L)$ to $G(K/k)$. 

It follows that the injection of $G(KL/K)$ into $G(L/k)$ by $\pi_2$ is also an isomorphism and that $\pi$
is an isomophism when $p$ is odd. As in the proof of theorem 1 it follows that the restriction to $H$ of the inflation of the cohomology class defining $D(p)$ must be zero, which is a contradiction, proving theorem 2
when $p$ is odd.

\bibliographystyle{plain}
\bibliography{references}

\end{document}